\documentclass{elsart}

\usepackage{amssymb}

\begin{document}

\begin{frontmatter}



\title{Optimal quadrature formulas with derivatives in Sobolev
space}


\author{Kh.M. Shadimetov, A.R. Hayotov, F.A. Nuraliev}

\address{Institute of Mathematics, National University of Uzbekistan, 29, Do`rmon yo`li str.,
100125, Tashkent, Uzbekistan.}

\ead{hayotov@mail.ru, abdullo\_hayotov@mail.ru}

\begin{abstract}
In the present paper the problem of construction of optimal
quadrature formulas in the sense of Sard in the space
$L_2^{(m)}(0,1)$is considered. Here the quadrature sum consists of
values of the integrand at nodes and values of the first and the third
derivatives of the integrand at the end points of the integration
interval. The coefficients of optimal quadrature formulas are
found and  the norm of the optimal error functional is calculated
for arbitrary natural number $N$ and for any $m\geq 4$ using S.L.
Sobolev method which is based on discrete analogue of the
differential operator $d^{2m}/dx^{2m}$. In particular, for $m=4,\
5$  optimality of the classical Euler-Maclaurin quadrature formula
is obtained. Starting from $m=6$ new optimal quadrature formulas
are obtained.
\end{abstract}

\end{frontmatter}

\section{Introduction. Statement of the problem}

We consider the following general quadrature formula
$$
\int\limits_0^1 p(x) \varphi (x)dx \cong  \sum\limits_{\beta  =
0}^N\sum\limits_{j=0}^{\alpha} C_{\beta j}
\varphi^{(j)}(x_{\beta}) \eqno (1.1)
$$
with the error functional
 $$
\ell (x) = p(x)\varepsilon _{[0,1]} (x) - \sum\limits_{\beta  =
0}^N\sum\limits_{j=0}^{\alpha} (-1)^jC_{\beta j}  \delta ^{(j)}(x
- x_\beta ) \eqno (1.2)
$$
in a Banach space $B$. Here $C_{\beta j}$ are the coefficients and
$x_{\beta}$ are the nodes of the formula (1.1), $N=1,2,...$,
$\alpha=0,1,...$, $p(x)$ is a weight function, $\varepsilon
_{[0,1]} (x)$ is the characteristic function of the interval [0,1], $\delta (x)$
is the Dirac delta-function, $\varphi$ is an element of the
space $B$.

The difference
$$
\left( {\ell ,\varphi } \right)
=\int\limits_{-\infty}^{\infty}\ell(x)\varphi(x)dx=\int\limits_0^1
p(x) \varphi (x)dx -  \sum\limits_{\beta  =
0}^N\sum\limits_{j=0}^{\alpha} C_{\beta j}
\varphi^{(j)}(x_{\beta}) \eqno(1.3)
$$
is called   \emph{the error} of the quadrature formula (1.1).

By the Cauchy-Schwarz inequality
$$
|(\ell,\varphi)|\leq \|\varphi|B\|\cdot \|\ell|B^*\|
$$
the error (1.3) of the formula (1.1) is estimated with the help of
the norm of the error functional (1.2) in the conjugate space $B^*$,
i.e. by
    $$
\left\| \ell |B^*  \right\| = \mathop {\sup }\limits_{\left\|
{\varphi |B } \right\| = 1} |(\ell ,\varphi )|.
$$
Thus estimation of the error (1.3) of the quadrature formula (1.1)
on functions of the space $B$ is reduced to finding the norm of
the error functional $\ell$ in the conjugate space $B^*$.

Obviously  the norm of the error functional $\ell $ depends on
the coefficients and the nodes of the quadrature formula (1.1). The
problem of finding the minimum of  the norm of the error functional
$\ell$ by coefficients and by nodes is called
\emph{S.M. Nikol'skii problem}, and obtained formula is called
\emph{optimal quadrature formula in the sense of Nikol'skii}. This
problem was first considered by S.M. Nikol'skii \cite{SMNik50}, and continued by
many authors, see e.g. \cite{BlaCom,Boj05,CatCom,Chak,SMNikBook88,Zhens} and references therein. Minimization
of the norm of the error functional $\ell$ by coefficients when the
nodes are fixed is called \emph{Sard's problem}. And obtained
formula is called \emph{optimal quadrature formula in the sense of
Sard}. First this problem was investigated by A.Sard \cite{Sard}.

The results of this paper are related to Sard's problem. So here
we discuss some of the previous results about optimal quadrature
formulas in the sense of Sard which are closely connected to our
results.

There are several methods of construction of optimal quadrature
formulas in the sense of Sard such as spline method, $\varphi-$
function method (see e.g. \cite{BlaCom,Schoenb65}) and Sobolev's
method which is based on construction of discrete analogue of a
linear differential operator (see e.g. \cite{Sobolev06,SobVas}).
In the different spaces, based on these methods, the Sard's
problem was investigated by many authors, see, for example,
\cite{IBab,BlaCom,CatCom,Com72a,Com72b,GhOs,HayMilShad11,HayMilShad123,Koh,FLan,MalOrl,
MeySar,Sard63,Schoenb65,Schoenb66,SchSil,Shad83,Shad02,ShadHay11,ShadHay120,ShadHay121,ShadHay122,Sobolev74,Sobolev06,SobVas,Zag,ZhaSha}
and references therein.

In the paper \cite{Schoenb65}, using spline method, optimality of
the classical Euler-Maclaurin formula was proved and the error of
this quadrature formula is calculated in $L_2^{(m)}(0,n)$, where
$L_2^{(m)}(0,n)$ is the space of functions which are square
integrable with $m$-th generalized derivative.

Let $W_{L_p}^m$ ($m=1,2,...,$ $1\leq p\leq \infty$) be a class of
functions $f$, having on the [0,1] $(m-1)$- absolute continues
derivative and $\|f^{(m)}\|_p\leq 1$, where $\|\cdot \|_p=\| \cdot
\|_{L_p(0,1)}$. In \cite{Zhens} it is proved, that  among
quadrature formulas (1.1) when $p(x)=1$ the Euler-Maclaurin
quadrature formula is optimal in the space $W_{L_p}^m$. And in
\cite{Shad01} optimality of the lattice cubature formulas of
Euler-Maclaurin type is proved in the space $L_2^{(m)}$.

Using $\varphi$-function method optimality of the Euler-Maclaurin
quadrature formula is proved and the error of this formula is
calculated by T. Catina\c{s} and Gh. Coman \cite{CatCom} in the space
$L_2^{(2)}(0,1)$. Also using this method in \cite{FLan} a procedure of
construction of quadrature formulas of the form (1.1), which are
exact for solutions of linear differential equations and are
optimal in the sense of Sard is discussed.

It should be noted, that in applications the formula (1.1) is
interesting for small values of $\alpha$. Optimal quadrature
formulas in the sense of Sard for the case $\alpha=0$ has already
been discussed by many authors, mainly in the space $L_2^{(m)}$
(see \cite{IBab,BlaCom,CatCom,Com72a,Com72b,GhOs,HayMilShad11,Koh,FLan,MalOrl,
MeySar,Sard63,Schoenb65,Schoenb66,SchSil,Shad83,Shad02,ShadHay11,Sobolev74,Sobolev06,SobVas,Zag,ZhaSha} and references therein).

The main aim of this paper is to construct optimal quadrature
formulas of the form (1.1) in the sense of Sard for the case
$\alpha=3$ when $p(x)=1$ in the space $L_2^{(m)}(0,1)$ equipped
with the norm
$$
\|\varphi(x)\|_{L_2^{(m)}(0,1)}=
\left\{\int\limits_0^1(\varphi^{(m)}(x))^2dx\right\}^{1/2}
$$
and $\int\limits_0^1(\varphi^{(m)}(x))^2dx<\infty$.

We use the Sobolev method \cite{Sobolev06,SobVas} which is
based on the discrete analogue of the differential operator
$d^{2m}/dx^{2m}$. We consider the following quadrature formula
$$
\int\limits_0^1 {\varphi (x)dx \cong \sum\limits_{\beta  = 0}^N {C
[\beta ]\varphi [\beta ] + A\bigg(\varphi '(0)-\varphi
'(1)\bigg)+B\bigg(\varphi '''(0)-\varphi''' (1)\bigg)\,\,} } \eqno
(1.4)
$$
with the error functional
 $$
\ell (x) = \varepsilon _{[0,1]} (x) - \sum\limits_{\beta  = 0}^N
{C[\beta ]\delta (x - h\beta ) + A\bigg(\delta '(x)-\delta
'(x-1)\bigg)+ B\bigg(\delta'''(x)-\delta'''(x-1)\bigg)} \eqno
(1.5)
$$
in the space $L_2^{(m)}(0,1)$ for $m\geq 4$. Here $C[\beta ]$,
$\beta = \overline{0,N}$,  $A$ and $B$ are the coefficients of the
formula (1.4),  $h = \frac{1}{N},$ $N$ is a natural number.

For the error functional (1.5) to be defined on the space
$L_2^{(m)}(0,1)$ it is necessary to impose the following
conditions (see \cite{Sobolev74})
$$
(\ell(x),x^{\alpha})=0,\ \ \ \alpha=0,1,2,...,m-1. \eqno(1.6)
$$
Hence it is clear that for existence of the quadrature formulas of
the form (1.4) the condition $N\ge m-3$ has to be met.

Note that here in after $\ell$ means the functional (1.5).

As was noted above by the Cauchy-Schwarz inequality, the error of the
formula (1.4) is estimated by the norm $\|\ell|L_2^{(m)*}(0,1)\|$
of the error functional (1.5). Furthermore the norm of the error
functional (1.5) depends on the coefficients $C[\beta ],$ $A$ and
$B$. We minimize the norm of the error functional (1.5) by the
coefficients $C[\beta]$, $A$ and $B$, i.e., we   find
$$
\left\| {\mathop \ell \limits^ \circ  |L_2^{(m)*} } \right\|  =
\mathop {\inf }\limits_{C[\beta ],A,B} \left\| {\ell|L_2^{(m)*} }
\right\|. \eqno     (1.7)
$$

The coefficients $C[\beta]$, $A$ and $B$ which satisfy the
equality (1.7) is called \emph{the optimal coefficients} and
denoted by $\stackrel{\circ}{C}[\beta]$, $\stackrel{\circ}{A}$ and
$\stackrel{\circ}{B}$ and the corresponding quadrature formula is
called \emph{the optimal quadrature formula in the sense of Sard}.
In the sequel, for the purposes of convenience the optimal
coefficients $\stackrel{\circ}{C}[\beta]$, $\stackrel{\circ}{A}$
and $\stackrel{\circ}{B}$ will be denoted as $C[\beta]$, $A$ and
$B$.

Thus to construct optimal quadrature formulas in the form (1.4) in
the sense of Sard we have to consequently solve the following
problems.

\textbf{Problem 1.} \emph{Find the norm of the error functional (1.5)
of the quadrature formula of the form (1.4) in the space
$L_2^{(m)*} (0,1)$.}

\textbf{Problem 2.} \emph{Find coefficients $C[\beta ],\,\,\,A$
and $B$ which satisfy the equality (1.7).}

The paper is organized as follows. In section 2 we give some definitions and known formulas.
In section 3 we determine the extremal function which corresponds to the error functional
$\ell$ and give a representation of the norm of the error
functional (1.5). Section 4 is devoted to a
minimization of $\left\| \ell \right\|^2 $ with respect to the
coefficients $C[\beta]$, $A$ and $B$. We obtain a system of linear equations  for
the coefficients of the optimal quadrature formula of the form (1.4) in the sense of
Sard in the space $L_2^{(m)}(0,1)$. Explicit
formulas for  coefficients of the optimal quadrature formula of
the form (1.4) are found in subsection 5.1. Moreover we
calculate the norm of the error functional
(1.5) of the optimal quadrature formula of the form (1.4) in subsection 5.2.

\section{Definitions and known formulas}

In this section we give some definitions and formulas that we need
to prove the main results.

Here the main concept used is that of functions of discrete
arguments and operations on them (see. \cite{Sobolev74,SobVas}).
For the purposes of completeness we give some definitions about
functions of discrete argument.

Assume that $\varphi $ and $\psi $ are real-valued functions
of real variable and are defined in real line $\mathbb{R}$.

\textbf{Definition 2.1.} Function $\varphi (h\beta )$ is called
\emph{function of discrete argument}, if it is given on some set
of integer values of $\beta $.

\textbf{Definition 2.2.} \emph{The inner product} of two discrete
functions $\varphi (h\beta )$ and $\psi (h\beta )$ is called the
number
    $$
\left[ {\varphi ,\psi } \right] = \sum\limits_{\beta  =  - \infty
}^\infty  {\varphi (h\beta ) \cdot \psi (h\beta )},
$$
if the series on the right hand side of the last equality converges
absolutely.

\textbf{Definition 2.3.} \textit{The convolution} of two discrete
functions $\varphi (h\beta )$ and $\psi (h\beta )$ is called the
inner product
$$
\varphi (h\beta )*\psi (h\beta ) = \left[ {\varphi (h\gamma
),\psi (h\beta  - h\gamma )} \right] = \sum\limits_{\gamma  =  -
\infty }^\infty  {\varphi (h\gamma ) \cdot \psi (h\beta  - h\gamma
)}.
$$

\emph{The Euler-Frobenius polynomials} $E_k (x)$, $k = 1,2,...$ is
defined by the following formula \cite{SobVas}
$$
E_k (x) = \frac{{(1 - x)^{k + 2} }}{x}\left( {x\frac{d}{{dx}}}
\right)^k \frac{x}{{(1 - x)^2 }},\eqno     (2.1)
$$
$E_0 (x) = 1$.

For the Euler-Frobenius polynomials $E_k(x)$ the following identity holds
$$
E_k (x) = x^k E_k \left( {\frac{1}{x}} \right),  \eqno (2.2)
$$
and also the following theorem is true

\textbf{Theorem 2.1} {\it (Lemma 3 of \cite{Shad10}). {Polynomial
$Q_k (x)$ which is defined by the formula}
$$
Q_k (x) = (x - 1)^{k + 1} \sum\limits_{i = 0}^{k + 1}
{\frac{{\Delta ^i 0^{k + 1} }}{{(x - 1)^i }}}\eqno (2.3)
$$
{is the Euler-Frobenius polynomial (2.1) of degree $k$, i.e.
$Q_k(x) = E_k(x)$, where $\Delta^i0^k=\sum_{l=1}^i(-1)^{i-l}C_i^ll^k.$} }

The following formula is valid \cite{Ham62}:
$$
\sum\limits_{\gamma  = 0}^{n - 1} {q^\gamma  \gamma ^k  =
\frac{1}{{1 - q}}\sum\limits_{i = 0}^k {\left( {\frac{q}{{1 - q}}}
\right)^i \Delta ^i 0^k  - \frac{{q^n }}{{1 - q}}\sum\limits_{i =
0}^k {\left( {\frac{q}{{1 - q}}} \right)^i \Delta ^i \gamma ^k
|_{\gamma  = n} ,} } }\eqno (2.4)
$$
where $\Delta^i\gamma^k$ is the finite difference of order $i$ of
$\gamma^k$, $q$ is ratio of a geometric progression.\\
At last we give the following well known formulas from \cite{Gel67}
$$
\sum\limits_{\gamma  = 0}^{\beta  - 1} {\gamma ^k  =
\sum\limits_{j = 1}^{k + 1} {\frac{{k!\,B_{k + 1 - j} }}{{j!\,(k +
1 - j)!}}\,\beta ^j ,} }\eqno    (2.5)
$$
where $B_{k + 1 - j} $ are Bernoulli numbers,
    $$
\Delta ^\alpha  x^\nu   = \sum\limits_{p = 0}^\nu C_\nu^p\Delta
^\alpha 0^p x^{\nu  - p}. \eqno (2.6)
$$

\section{The extremal function and the representation of the error functional norm}

To solve Problem 1, i.e., for finding the norm of the error
functional (1.5) in the space $L_2^{(m)}(0,1)$ a concept of the
extremal function is used \cite{Sobolev74}. The function $\psi
_\ell $ is said to be \emph{the extremal function} of the error
functional (1.5) if the following equality holds
$$
\left( {\ell,\psi _\ell} \right) = \left\| {\ell
|L_2^{(m)*}} \right\|\left\| {\psi _\ell |L_2^{(m)} }
\right\|. \eqno (3.1)
$$
In the space $L_2^{(m)}$ the extremal function $\psi _\ell$ of
a functional $\ell $ was found by S.L. Sobolev \cite{Sobolev74,SobVas}. This
extremal function has the form
$$
\psi _\ell(x)= ( - 1)^m \ell (x) * G(x) + P_{m - 1} (x),\eqno
(3.2)
$$
where
$$
G(x)=\frac{{|x|^{2m - 1}}}{{2\cdot(2m - 1)!}}\eqno (3.3)
$$
is a solution of the equation
$$
\frac{d^{2m}}{dx^{2m}}G(x)=\delta(x),\eqno (3.4)
$$
$P_{m - 1}(x)$ is a polynomial of degree $m - 1$, the symbol
*  is operation of convolution, i.e.
$$
f(x)*g(x) = \int\limits_{ - \infty }^\infty  {f(x - y)g(y)dy = }
\int\limits_{ - \infty }^\infty  {f(y)g(x - y)dy}.
$$

It is well known that for any functional $\ell$ in $L_2^{(m)*}$
the equality
\begin{eqnarray*}
\|\ell|L_2^{(m)*}(0,1)\|^2&=&(\ell,\psi_{\ell})=\left(\ell(x),(-1)^m\ell(x)*G(x)\right)=\\
&=&\int\limits_{-\infty}^{\infty}\ell(x)
\left((-1)^m\int\limits_{-\infty}^{\infty}\ell(y)G(x-y)dy\right)dx
\end{eqnarray*}
holds \cite{Sobolev74}.\\
Applying this equality to the error functional (1.5) we obtain the
following
$$
\left\| \ell  \right\|^2=(\ell ,\psi _\ell  ) =$$ $$=( - 1)^m\left[
{\sum\limits_{\beta  = 0}^N {\sum\limits_{\gamma  = 0}^N {C[\beta
]C[\gamma ]} } } \right.{{|h\beta  - h\gamma |^{2m -1} } \over
{2(2m - 1)!}} -2\sum\limits_{\beta  = 0}^N C[\beta ] \int\limits_0^1 {\frac{{\left| {x - h\beta }
\right|^{2m - 1} }}{{2(2m - 1)!}}} dx
-
$$
$$
-2A\sum\limits_{\beta=0}^NC[\beta]P_{2m-2}[\beta]-2B\sum\limits_{\beta=0}^NC[\beta]P_{2m-4}[\beta]+\frac{2A}{(2m-1)!}+
$$
$$
+\frac{2B+A^2}{(2m-3)!}+\frac{2AB}{(2m-5)!}+\frac{B^2}{(2m-7)!} -
\left. \frac{1}{(2m+1)!} \right]. \eqno (3.5)
$$
where $P_{k}(h\beta)=\frac{(h\beta)^k+(1-h\beta)^k}{2k!}$.

 Thus Problem 1 is solved for quadrature formulas of the form
(1.4) in the space $L_2^{(m)}(0,1)$.

\section{The  system  for optimal coefficients of the quadrature formula (1.4)}

Now we investigate  Problem 2. For finding the minimum of
$\left\| \ell \right\|^2 $ under the conditions (1.6) the Lagrange
method is used. For this we consider the following function
    $$
\Psi  = \left\|\ell\right\|^2  - 2 \cdot ( - 1)^m
\sum\limits_{\alpha  = 0}^{m - 1} {\lambda _\alpha  }
(\ell(x),x^\alpha  ),
$$
where $\lambda_{\alpha}$ are unknown multipliers. The function
$\Psi$ is the multidimensional function with respect to the
coefficients $C[\beta],\ A,\ B$ and $\lambda_{\alpha}$. Equating
to zero partial derivatives of the $\Psi$ by coefficients $C[\beta
],\,\,A$ and $B,$ together with the conditions (1.6) we get the
following system of linear equations
                $$
\sum\limits_{\gamma  = 0}^N {C[\gamma ]} \frac{{\left| {h\beta  - h\gamma } \right|^{2m - 1} }}{{2(2m - 1)!}} -
 AP_{2m-2}(h\beta)- $$ $$-BP_{2m-4}(h\beta)+\sum\limits_{\alpha =
0}^{m - 1}
 \,\lambda _\alpha(h\beta)^{\alpha}=f_{m}(h\beta),\,\,\,\,\,\,\beta  =
\overline {0,N} , \eqno   (4.1)
$$
$$
\sum\limits_{\beta=0}^NC[\beta]P_{2m-2}(h\beta)-\frac{A}{(2m-3)!}-\frac{B}{(2m-5)!}
+\sum\limits_{\alpha=2}^{m-1}\alpha\cdot\lambda_{\alpha}=\frac{1}{(2m-1)!},\eqno
(4.2)
$$
$$
\sum\limits_{\beta=0}^NC[\beta]P_{2m-4}(h\beta)-\frac{A}{(2m-5)!}-\frac{B}{(2m-7)!}+$$
$$+\sum\limits_{\alpha=4}^{m-1}
\alpha(\alpha-1)(\alpha-2)\cdot\lambda_{\alpha}=\frac{1}{(2m-3)!},\eqno
(4.3)
$$
        $$
\sum\limits_{\beta=0}^NC[\beta](h\beta)^i=\frac{1}{i+1},\ \ \
i=\overline{0,1}\eqno (4.4)
$$
                 $$
\sum\limits_{\beta=0}^NC[\beta](h\beta)^j-jA=\frac{1}{j+1},\ \ \
j=\overline{2,3}\eqno  (4.5)
$$
$$
\sum\limits_{\beta=0}^NC[\beta](h\beta)^{\alpha}-\alpha
A-\alpha(\alpha-1)(\alpha-2)B=\frac{1}{\alpha+1},\ \ \
\alpha=\overline{4,m-1}.\eqno     (4.6)
$$
where $P_{k}(h\beta)=\frac{(h\beta)^k+(1-h\beta)^k}{2k!}$.

 The
system (4.1)-(4.6) is called   \emph{the discrete system of
Wiener-Hopf type} for the optimal coefficients
\cite{Sobolev74,SobVas}. In the system (4.1)-(4.6) the
coefficients $C[\beta ],$ $\beta = \overline {0,N}$, $A$ and $B,$
and also $\lambda _\alpha ,\,\,\,\,\alpha  = \overline {0,m - 1}$
are unknowns. The system (4.1)-(4.6) has unique solution and this
solution gives the minimum to the $\|\ell\|^2$. Here we omitted
the proof of the existence and uniqueness of the solution of the
system (4.1) - (4.6). The proof of the existence and uniqueness of
the solution of this system is as the proof of the existence and
uniqueness of the solution of discrete Wiener-Hopf type system of
the optimal coefficients in the space $L_2^{(m)}(0,1)$ for
quadrature formulas of the form (1.1) for the case $\alpha=0$ (see
\cite{Sobolev74,SobVas}). It should be noted, that in \cite{FLan}
the uniqueness of the optimal quadrature formulas in the Sard's
sense of the form (1.1) is discussed.

\section{The coefficients and the norm of the error functional of
the optimal quadrature formulas}

In the present section we study the solution of the system
(4.1)-(4.6). To solve this system we use the approach which was
suggested by S.L. Sobolev in \cite{Sobolev06}.

\subsection{The coefficients of the optimal quadrature formulas of the form (1.4)}

Suppose that $C[\beta]=0$  for $\beta  < 0$ and $\beta  > N$.
Using Definition 2.3  we rewrite the equation (4.1) in the
convolution form:
                $$
C[\beta ]*\frac{{\left| h\beta \right|^{2m - 1} }}{{2(2m - 1)!}} -
 AP_{2m-2}(h\beta)- BP_{2m-4}(h\beta)+\sum\limits_{\alpha =
0}^{m - 1}
 \,\lambda _\alpha(h\beta)^{\alpha}=f_{m}(h\beta), \,\,\,\,\,\,\beta=\overline{0,N}, \eqno
(5.1)
$$
where
$$
f_{m}(h\beta)= \int\limits_0^1 {\frac{{\left| {x - h\beta }
\right|^{2m - 1} }}{{2(2m - 1)!}}} dx=\frac{(h\beta)^{2m}}{(2m)!}+
\sum\limits_{j=0}^{2m-1}\frac{(-h\beta)^{2m-1-j}}{2(2m-1-j)!\cdot
(j+1)!}. \eqno (5.2)
$$

We consider the following problem.

\textbf{Problem A.} \emph{Find the discrete function $C[\beta]$
and unknown coefficients $A,$ $B$, $\lambda_{\alpha}$, which
satisfy the system (4.1)-(4.6).}

Further, instead of $C[\beta]$ we introduce the functions
$$
v(h\beta) =C[\beta ]*\frac{{\left| h\beta \right|^{2m - 1}
}}{{2(2m - 1)!}},\eqno     (5.3)
$$
$$
u(h\beta) = v(h\beta) - AP_{2m-2}(h\beta)-
BP_{2m-4}(h\beta)+\sum\limits_{\alpha = 0}^{m - 1}
 \,\lambda _\alpha(h\beta)^{\alpha}. \eqno  (5.4)
$$

In this statement it is necessary to express $C[\beta]$ by the
function $u(h\beta)$. For this we need such operator
$D_m(h\beta)$, which satisfies the equation
$$
hD_m(h\beta)*G(h\beta) = \delta(h\beta),\eqno    (5.5)
$$
where $G(h\beta)=\frac{{\left| h\beta \right|^{2m - 1} }}{{2(2m -
1)!}}$ is the discrete argument function corresponding to $G(x)$
defined by (3.3), $\delta(h\beta)$ is equal to 0 when $\beta \ne
0$ and is equal to 1 when $\beta = 0$, i.e. $\delta(h\beta)$ is
the discrete delta-function. The equation (5.5) is the discrete
analogue of equation (3.4). So the discrete function
$D_m(h\beta)$ is called \emph{the discrete analogue} of the
differential operator $d^{2m}/dx^{2m}$ \cite{Sobolev74}.

It should be noted that the operator $D_m(h\beta)$ was firstly
introduced and investigated by S.L. Sobolev \cite{Sobolev74}.

In \cite{Shad85} the discrete analogue $D_m(h\beta)$ of the
differential operator $d^{2m}/dx^{2m}$, which satisfies  equation
(5.5), is constructed and the following theorem is proved.

\textbf{Theorem 5.1.} {\it The discrete analogue of the
differential operator ${d^{2m}}/{dx^{2m}}$ has the form
$$
D_m(h\beta)=\frac{(2m-1)!}{h^{2m}}\left\{
\begin{array}{lll}
{\displaystyle
\sum\limits_{k=1}^{m-1}\frac{(1-q_k)^{2m+1}q_k^{|\beta|}}
{q_kE_{2m-1}(q_k)}        }
& \mbox{ for }& |\beta|\geq 2,\\
{\displaystyle 1+\sum\limits_{k=1}^{m-1}\frac{(1-q_k)^{2m+1}}
{E_{2m-1}(q_k)}                  }
& \mbox{ for }& |\beta|= 1,\\
{\displaystyle
-2^{2m-1}+\sum\limits_{k=1}^{m-1}\frac{(1-q_k)^{2m+1}} {q_k
E_{2m-1}(q_k)}           } & \mbox{ for }& \beta= 0,
\end{array}
\right. \eqno (5.6)
$$
where $E_{2m-1}(q)$ is the Euler-Frobenius polynomial of degree
$2m-1$, $q_k$ are the roots of the Euler-Frobenius polynomial
$E_{2m-2}(q)$, $|q_k|<1$, $h$ is a small positive parameter.}

Furthermore several properties  of the discrete argument function $D_m(h\beta)$ were proved in \cite{Shad85}.
Here we give the following property of the discrete argument function $D_m(h\beta)$ which we need in our
computations.

{\bf Theorem 5.2.} {\it The discrete argument function
$D_m(h\beta)$ and the monomials $(h\beta)^k$ are related to each
other as follows}
$$
\sum_{\beta=-\infty}^{\infty}D_m(h\beta)(h\beta)^k= \left\{
\begin{array}{lll}
0&\mbox{\it when } & 0\leq k\leq 2m-1,\\
(2m)!&\mbox{\it when } & k= 2m,
\end{array}
\right. \eqno (5.7)
$$
$$
\sum_{\beta=-\infty}^{\infty}D_m(h\beta)(h\beta)^k= \left\{
\begin{array}{lll}
0&\mbox{\it when } & 2m+1\leq k\leq 4m-1,\\
{\displaystyle             \frac{h^{2m}(4m)!B_{2m}}{(2m)!}    }
&\mbox{\it when }  & k= 4m.
\end{array}
\right.
$$

Then, taking into account (5.5) and Theorems 5.1, 5.2, for the
optimal coefficients $C[\beta]$ we have
$$
C[\beta]=hD_m(h\beta)*u(h\beta).\eqno     (5.8)
$$

So, if we find the function $u(h\beta)$, then the optimal
coefficients $C[\beta]$ will be found from  equality (5.8).

To calculate the convolution (5.8) it is required to find the
representation of the function $u(h\beta)$ for all integer values
of $\beta$. From the equality (5.1) we get, that $u(h\beta)
=f_m(h\beta)$ when $h\beta \in [0,1]$, where $f_m(h\beta)$ is defined
by equality (5.2). Now we need to find the representation of the
function $u(h\beta)$ when $\beta < 0$ and $\beta>N$.

Since $C[\beta]= 0$ when $h\beta  \notin [0,1]$, then
$$
C[\beta]=hD_m(h\beta)*u(h\beta) = 0,\,\,\,\,\,\,\,\,\,h\beta
\notin [0,1].
$$

Now we calculate the convolution
$v(h\beta)=C[\beta]*\frac{|h\beta|^{2m-1}}{2(2m-1)!}$ when
$h\beta\notin [0,1]$.

Suppose $\beta<0$, then taking into account (4.4)-(4.6), we have
$$
v(h\beta)=C[\beta]*\frac{|h\beta|^{2m-1}}{2(2m-1)!}=\sum\limits_{\gamma=-\infty}^{\infty}
C[\gamma]\frac{|h\beta-h\gamma|^{2m-1}}{2(2m-1)!}=
$$
$$
=-\sum\limits_{\alpha=0}^{m-1}\frac{(h\beta)^{2m-1-\alpha}(-1)^{\alpha}}
{2(2m-1-\alpha)!\alpha!}\sum\limits_{\gamma=0}^NC[\gamma](h\gamma)^{\alpha}-
\sum\limits_{\alpha=m}^{2m-1}\frac{(h\beta)^{2m-1-\alpha}(-1)^{\alpha}}
{2(2m-1-\alpha)!\alpha!}\sum\limits_{\gamma=0}^NC[\gamma](h\gamma)^{\alpha}=
$$
$$
=-\frac{(h\beta)^{2m-1}}{2(2m-1)!}+\frac{(h\beta)^{2m-2}}{4(2m-2)!}-\frac{(h\beta)^{2m-3}}{4(2m-3)!}\left(\frac{1}{3}+2A\right)+
\frac{(h\beta)^{2m-4}}{2(2m-4)!3!}\left(\frac{1}{4}+3A\right)-$$
$$-
\sum\limits_{\alpha=4}^{m-1}\frac{(h\beta)^{2m-1-\alpha}(-1)^{\alpha}}{2(2m-1-\alpha)\alpha!}
\left(\frac{1}{\alpha+1}+\alpha
 A+\alpha(\alpha-1)(\alpha-2)B\right)-
$$
$$
-\sum\limits_{\alpha=m}^{2m-1}\frac{(h\beta)^{2m-1-\alpha}(-1)^{\alpha}}
{2(2m-1-\alpha)!\alpha!}\sum\limits_{\gamma=0}^NC[\gamma](h\gamma)^{\alpha}.
$$
Hence, denoting by
$R_{m-1}(h\beta)=\sum\limits_{\alpha=m}^{2m-1}\frac{(h\beta)^{2m-1-\alpha}(-1)^{\alpha}}
{2(2m-1-\alpha)!\alpha!}\sum\limits_{\gamma=0}^NC[\gamma](h\gamma)^{\alpha}$
for the case  $\beta<0$ we get
$$
v(h\beta)=
-\frac{(h\beta)^{2m-1}}{2(2m-1)!}+\frac{(h\beta)^{2m-2}}{4(2m-2)!}-\frac{(h\beta)^{2m-3}}{4(2m-3)!}\left(\frac{1}{3}+2A\right)+
\frac{(h\beta)^{2m-4}}{2(2m-4)!3!}\left(\frac{1}{4}+3A\right)-
$$
$$-\sum\limits_{\alpha=4}^{m-1}\frac{(h\beta)^{2m-1-\alpha}(-1)^{\alpha}}{2(2m-1-\alpha)\alpha!}
\left(\frac{1}{\alpha+1}+\alpha
 A+\alpha(\alpha-1)(\alpha-2)B\right)-R_{m-1}(h\beta). \eqno
(5.9)
$$
Now suppose $\beta>N$ then for $v(h\beta)$ we get
$$
v(h\beta)=
\frac{(h\beta)^{2m-1}}{2(2m-1)!}-\frac{(h\beta)^{2m-2}}{4(2m-2)!}+\frac{(h\beta)^{2m-3}}{4(2m-3)!}\left(\frac{1}{3}+2A\right)-
\frac{(h\beta)^{2m-4}}{2(2m-4)!3!}\left(\frac{1}{4}+3A\right)+
$$
$$+
\sum\limits_{\alpha=4}^{m-1}\frac{(h\beta)^{2m-1-\alpha}(-1)^{\alpha}}{2(2m-1-\alpha)\alpha!}
\left(\frac{1}{\alpha+1}+\alpha
 A+\alpha(\alpha-1)(\alpha-2)B\right)+R_{m-1}(h\beta). \eqno
(5.10)
$$
Denoting
$$
R_{m-1}^{(-)}(h\beta)=\sum\limits_{\alpha=0}^{m-1}
\lambda_{\alpha}(h\beta)^{\alpha}-R_{m-1}(h\beta), \eqno (5.11)
$$
$$
R_{m-1}^{(+)}(h\beta)=\sum\limits_{\alpha=0}^{m-1}
\lambda_{\alpha}(h\beta)^{\alpha}+R_{m-1}(h\beta), \eqno (5.12)
$$
and taking into account (5.9), (5.10), (5.4) we have the following
problem

\textbf{Problem B.} \emph{Find the solution of the equation}
$$
hD_m(h\beta)*u(h\beta) = 0,\,\,\,\,\,\,\,\,h\beta  \notin
[0,1]\eqno (5.13)
$$
\emph{having the form:}

$$
u(h\beta) = \left\{
\begin{array}{ll}
-\frac{(h\beta)^{2m-1}}{2(2m-1)!}+\frac{(h\beta)^{2m-2}}{4(2m-2)!}-\frac{(h\beta)^{2m-3}}{4(2m-3)!}\left(\frac{1}{3}+2A\right)+
\frac{(h\beta)^{2m-4}}{2(2m-4)!3!}\left(\frac{1}{4}+3A\right)- & \\
-\sum\limits_{\alpha=4}^{m-1}\frac{(h\beta)^{2m-1-\alpha}(-1)^{\alpha}}{2(2m-1-\alpha)\alpha!}
\left(\frac{1}{\alpha+1}+\alpha
 A+\alpha(\alpha-1)(\alpha-2)B\right)+R_{m-1}^{(-)}(h\beta),& \beta<0,\\
f(h\beta), &0\leq \beta\leq N,\\
\frac{(h\beta)^{2m-1}}{2(2m-1)!}-\frac{(h\beta)^{2m-2}}{4(2m-2)!}+\frac{(h\beta)^{2m-3}}{4(2m-3)!}\left(\frac{1}{3}+2A\right)-
\frac{(h\beta)^{2m-4}}{2(2m-4)!3!}\left(\frac{1}{4}+3A\right)+ & \\
+
\sum\limits_{\alpha=4}^{m-1}\frac{(h\beta)^{2m-1-\alpha}(-1)^{\alpha}}{2(2m-1-\alpha)\alpha!}
\left(\frac{1}{\alpha+1}+\alpha
 A+\alpha(\alpha-1)(\alpha-2)B\right)+R_{m-1}^{(+)}(h\beta),& \beta>N,\\
\end{array}
\right.\eqno (5.14)
$$

\emph{where $R_{m-1}^{(-)}(h\beta)$ and $R_{m-1}^{(+)}(h\beta)$
are unknown polynomials of degree $m-1$ and $A$, $B$ are unknown
coefficients.}\\ If we find $R_{m-1}^{(-)}(h\beta)$ and
$R_{m-1}^{(+)}(h\beta)$, then from (5.11), (5.12) we obtain
    $$
\sum\limits_{\alpha=0}^{m-1}\lambda_{\alpha}(h\beta)^{\alpha}=\frac{1}{2}
\left(R_{m-1}^{(+)}(h\beta)+R_{m-1}^{(-)}(h\beta)\right),
$$
  $$
R_{m-1}(h\beta)=\frac{1}{2}
\left(R_{m-1}^{(+)}(h\beta)-R_{m-1}^{(-)}(h\beta)\right).
$$

Unknowns $R_{m-1}^{(-)}(h\beta)$, $R_{m-1}^{(+)}(h\beta)$, $A$ and
$B$ can be found from equation (5.13), using the discrete argument
function $D_m(h\beta)$. Then we can obtain the explicit form of
the function $u(h\beta)$ and respectively  we can find the optimal
coefficients $C[\beta]$ ($\beta=0,1,...,N$), $A$ and $B$. Thus
Problem B and respectively Problem A can be solved.

But here we will not find $R_{m-1}^{(-)}(h\beta)$,
$R_{m-1}^{(+)}(h\beta)$. Instead, using $D_m(h\beta)$ and the form
(5.14) of the discrete argument  function $u(h\beta)$, taking into
account (5.8), we find the expressions for the optimal
coefficients $C[\beta]$ when $\beta=1,2,...,N-1$.

We introduce the following notations
$$
d_k=\frac{(2m-1)!(1-q_k)^{2m+1}}{h^{2m}q_kE_{2m-1}(q_k)}
\sum\limits_{\gamma=1}^{\infty}q_k^{\gamma}\Bigg\{\frac{(h\gamma)^{2m-1}}{2(2m-1)!}+
\frac{(h\gamma)^{2m-2}}{4(2m-2)!}-
\frac{(h\gamma)^{2m-3}}{4(2m-3)!}\left(\frac{1}{3}+2A\right)+
$$
$$
+\frac{(h\gamma)^{2m-4}}{12(2m-4)!}\left(\frac{1}{4}+3A\right)+
\sum\limits_{\alpha=4}^{m-1}\frac{(h\gamma)^{2m-1-\alpha}(-1)^{\alpha}}{2(2m-1-\alpha)\alpha!}
\left(\frac{1}{\alpha+1}+\alpha
 A+\alpha(\alpha-1)(\alpha-2)B\right) +
$$
$$
+R_{m-1}^{(-)}(-h\gamma)-f(-h\gamma)\Bigg\}, \eqno (5.15)
$$
$$
p_k=\frac{(2m-1)!(1-q_k)^{2m+1}}{h^{2m}q_kE_{2m-1}(q_k)}
\sum\limits_{\gamma=1}^{\infty}q_k^{\gamma}\Bigg\{\frac{(h(\gamma+N))^{2m-1}}{2(2m-1)!}+
\frac{(h(\gamma+N))^{2m-2}}{4(2m-2)!}-$$
$$
-\frac{(h(\gamma+N))^{2m-3}}{4(2m-3)!}\left(\frac{1}{3}+2A\right)
+\frac{(h(\gamma+N))^{2m-4}}{12(2m-4)!}\left(\frac{1}{4}+3A\right)+$$ $$+
\sum\limits_{\alpha=4}^{m-1}\frac{(h(\gamma+N))^{2m-1-\alpha}(-1)^{\alpha}}{2(2m-1-\alpha)\alpha!}
\left(\frac{1}{\alpha+1}+\alpha
 A+\alpha(\alpha-1)(\alpha-2)B\right) +
$$
$$
+R_{m-1}^{(+)}(h\gamma+1)-f(h\gamma+1)\Bigg\}, \eqno (5.16)
$$
where $k=1,2,...,m-1$, $E_{2m-1}(q)$ is the Euler-Frobenius
polynomial of degree $2m-1$, $q_k$ are given in Theorem 5.1. Note
that because of $|q_k|<1$ the series in the (5.15) and (5.16) are
convergent.

The following holds

\textbf{Theorem 5.3.} \emph{The coefficients $C[\beta]$,
$\beta=1,2,...,N-1$ of the optimal quadrature formulas of the form
(1.4) in the space $L_2^{(m)}(0,1)$ for $m\geq 4$ have the following form}
$$
C[\beta]=h\left(1+\sum\limits_{k=1}^{m-1}\left(d_kq_k^{\beta}+p_kq_k^{N-\beta}
\right)\right),\ \ \beta=1,2,...,N-1, \eqno (5.17)
$$
\emph{where $d_k, p_k$ are defined by (5.15), (5.16), $q_k$ are
given in Theorem 5.1.}

\textbf{Proof}.  Suppose $\beta=1,2,...,N-1$. Then from (5.8),
using Definition 2.3,  equalities (5.6), (5.14), we have
$$
C[\beta]=hD_m(h\beta)*u(h\beta)=h\sum\limits_{\gamma=-\infty}^{\infty}
D_m(h\beta-h\gamma)u(h\gamma)=
$$
$$
=h\left(\sum\limits_{\gamma=-\infty}^{-1}D_m(h\beta-h\gamma)u(h\gamma)+
\sum\limits_{\gamma=0}^ND_m(h\beta-h\gamma)f_m(h\gamma)+\sum\limits_{\gamma=N+1}^{\infty}
D_m(h\beta-h\gamma)u(h\gamma)\right).
$$
Now, adding and subtracting the expressions
$h\sum\limits_{\gamma=-\infty}^{-1}D_m(h\beta-h\gamma)f_m(h\gamma)$
and
$h\sum\limits_{\gamma=N+1}^{\infty}D_m(h\beta-h\gamma)f_m(h\gamma)$
to and from the last expression and taking into account Definition
2.3 we get
$$
C[\beta]=h\Bigg\{D_m(h\beta)*f_m(h\beta)+\sum\limits_{k=1}^{m-1}q_k^{\beta}\frac{(2m-1)!}{h^{2m}}
\frac{(1-q_k)^{2m+1}}{q_kE_{2m-1}(q_k)}\sum\limits_{\gamma=1}^{\infty}q_k^{\gamma}\Bigg[-
\frac{(h\gamma)^{2m-1}}{2(2m-1)!}+
$$
$$
+
\frac{(h\gamma)^{2m-2}}{4(2m-2)!}-
\frac{(h\gamma)^{2m-3}}{4(2m-3)!}\left(\frac{1}{3}+2A\right)+\frac{(h\gamma)^{2m-4}}{12(2m-4)!}\left(\frac{1}{4}+3A\right)-
$$
$$
-
\sum\limits_{\alpha=4}^{m-1}\frac{(h\gamma)^{2m-1-\alpha}(-1)^{\alpha}}{2(2m-1-\alpha)\alpha!}
\left(\frac{1}{\alpha+1}+\alpha
 A+\alpha(\alpha-1)(\alpha-2)B\right) +
$$
$$
+R_{m-1}^{(-)}(-h\gamma)-f_m(-h\gamma)\Bigg]+
$$
$$
+\sum\limits_{k=1}^{m-1}q_k^{N-\beta}\frac{(2m-1)!}{h^{2m}}
\frac{(1-q_k)^{2m+1}}{q_kE_{2m-1}(q_k)}\sum\limits_{\gamma=1}^{\infty}q_k^{\gamma}\Bigg[
\frac{(h(\gamma+N))^{2m-1}}{2(2m-1)!}-
$$
$$
-
\frac{(h(\gamma+N))^{2m-2}}{4(2m-2)!}
+\frac{(h(\gamma+N))^{2m-3}}{4(2m-3)!}\left(\frac{1}{3}+2A\right)
-\frac{(h(\gamma+N))^{2m-4}}{12(2m-4)!}\left(\frac{1}{4}+3A\right)+$$ $$+
\sum\limits_{\alpha=4}^{m-1}\frac{(h(\gamma+N))^{2m-1-\alpha}(-1)^{\alpha}}{2(2m-1-\alpha)\alpha!}
\left(\frac{1}{\alpha+1}+\alpha
 A+\alpha(\alpha-1)(\alpha-2)B\right) +
$$
$$
+R_{m-1}^{(+)}(h\gamma+1)-f_m(h\gamma+1)\Bigg]\Bigg\}.
$$
Hence taking into account the notations (5.15), (5.16) we obtain
$$
C[\beta]=h\left(D_m(h\beta)*f_m(h\beta)+\sum\limits_{k=1}^{m-1}(d_kq_k^{\beta}+
p_kq_k^{N-\beta})\right).\eqno (5.18)
$$
Now using Theorems 5.1, 5.2 and  equality (5.2) we get
$$
D_m(h\beta)*f(h\beta)=D_m(h\beta)*\left(\frac{(h\beta)^{2m}}{(2m)!}+
\sum\limits_{j=0}^{2m-1}\frac{(-h\beta)^{2m-1-j}}{2(2m-1-j)!(j+1)!}\right)
$$
$$
=D_m(h\beta)*\frac{(h\beta)^{2m}}{(2m)!}=1. \eqno (5.19)
$$
Putting (5.19) to the equation (5.18) we get (5.17). Theorem 5.3
is proved.

Furthermore we need the following lemmas in the proof of the main
results.

\textbf{Lemma 5.1.} \textit{The following identity is taken place
    $$
\sum\limits_{i = 0}^\alpha  {\frac{{d_kq_k + p_kq_k^{N + i} ( -
1)^{i + 1} }}{{(q_k - 1)^{i + 1} }}\Delta ^i 0^\alpha  }  = ( -
1)^{\alpha + 1} \sum\limits_{i = 0}^\alpha  {\frac{{d_kq_k^i  +
p_kq_k^{N + 1}(-1)^{i + 1} }}{{(1 - q_k)^{i + 1} }}\Delta ^i
0^\alpha  },\eqno (5.20)
$$
here $\alpha $ and $N$ are natural numbers, $d_k$ and $p_k$ are
defined by (5.15), (5.16), $\Delta ^i 0^\alpha$ is given in
Theorem 2.1, $q_k$ are given in Theorem 5.1.}

\textbf{Proof.}  For the purposes of convenience the left and the
right hand sides of (5.20) we denote by $L_1 $ and
$(-1)^{\alpha+1}L_2$ respectively, i.e.
    $$
L_1  = \sum\limits_{i = 0}^\alpha  {\frac{{d_kq_k + p_kq_k^{N + i}
( - 1)^{i + 1} }}{{(q_k  - 1)^{i + 1} }}\Delta ^i 0^\alpha  }
\mbox{ and }L_2  = \sum\limits_{i = 0}^\alpha {\frac{{d_kq_k^i  +
p_kq_k^{N + 1} (-1)^{i + 1} }}{{(1- q_k)^{i + 1} }}\Delta ^i
0^\alpha} .$$

First consider $L_1 $. Using the equality (2.3) and the identity
(2.2) for $L_1 $ consequently we get
$$
L_1  = \sum\limits_{i = 0}^\alpha  {\frac{{d_kq_k+ p_kq_k^{N + i}
( - 1)^{i + 1} }}{{(q_k  - 1)^{i + 1} }}\Delta ^i 0^\alpha  }  =
\frac{{d_kq_k}}{{(q_k - 1)^{\alpha  + 1} }}E_{\alpha  - 1}(q_k) +
\frac{{p_kq_k^{N + \alpha } ( - 1)^{\alpha  + 1} }}{{(q_k-
1)^{\alpha  + 1} }}E_{\alpha  - 1} \left( {\frac{1}{q_k}} \right)
=
$$
$$
= \frac{{d_kq_k}}{{(q_k - 1)^{\alpha  + 1} }}E_{\alpha  - 1}(q_k)
+ \frac{{p_kq_k^{N + \alpha}(-1)^{\alpha  + 1}}}{{(q_k -
1)^{\alpha + 1} }}\frac{{E_{\alpha  - 1}(q_k)}}{{q_k^{\alpha  - 1}
}} = \frac{{d_kq_k + p_kq_k^{N + 1}(- 1)^{\alpha  + 1} }}{{(q_k -
1)^{\alpha + 1} }}E_{\alpha  - 1}(q_k).\eqno (5.21)
$$
Similarly for $L_2 $   using  (2.3) and (2.2) we have
$$
L_2  = \sum\limits_{i = 0}^\alpha  {\frac{{d_kq_k^i  + p_kq_k^{N +
1}( - 1)^{i + 1} }}{{(1 - q_k)^{i + 1} }}\Delta ^i 0^\alpha  }  =
\frac{{d_kq_k^\alpha  }}{{(q_k - 1)^{\alpha  + 1} }}E_{\alpha  -
1} \left( {\frac{1}{q_k}} \right) + \frac{{p_kq_k^{N + 1} }}{{(q_k
- 1)^{\alpha  + 1} }}E_{\alpha  - 1} (q_k) =
$$
$$
= \frac{{d_kq_k}}{{(1 - q_k)^{\alpha  + 1} }}E_{\alpha  - 1}(q_k)
+ \frac{{p_kq_k^{N + 1} }}{{(q_k - 1)^{\alpha  + 1} }}E_{\alpha  -
1}(q_k)= \frac{{d_kq_k( - 1)^{\alpha  + 1}  + p_kq_k^{N + 1}
}}{{(q_k - 1)^{\alpha  + 1} }}E_{\alpha  - 1} (q_k) =
$$
$$
= ( - 1)^{\alpha  + 1} \frac{{d_kq_k + p_kq_k^{N + 1} ( -
1)^{\alpha + 1} }}{{(q_k-1)^{\alpha  + 1} }}E_{\alpha  - 1}
(q_k).\eqno (5.22)
$$
From (5.21) and (5.22) it is clear, that $L_1  = ( - 1)^{\alpha  +
1} L_2$. Lemma 5.1 is proved.

We denote
    $$
Z_p  = \sum\limits_{k = 1}^{m - 1} {\sum\limits_{i = 0}^p
{\frac{{d_k q_k^{N + i}  + p_k q_k ( - 1)^{i + 1} }}{{(1 - q_k
)^{i + 1} }}\Delta ^i 0^p } }. \eqno  (5.23)
$$

\textbf{Lemma 5.2.} \emph{ The following identities are valid}
\begin{eqnarray*}
\sum\limits_{j = 0}^{m - 1} {\frac{{( - 1)^{j + 1} }}{{(j + 1)!}}\sum\limits_{i = 1}^{2m - 2 - j} {\frac{{B_{2m
- j - i} h^{2m - j - i} }}{{i!\,\,(2m - j - i)!}} } }&=&  \sum\limits_{j = 2}^m {\frac{{B_j h^j
}}{{j!}}\sum\limits_{i = 1}^{m } {\frac{{( - 1)^i }}{{i!\,\,(2m + 1 - j - i)!}} + } } \\
&&+\sum\limits_{j = m + 1}^{2m - 1} {\frac{{B_j h^j }}{{j!}}\sum\limits_{i = 1}^{2m - j} {\frac{{( - 1)^i
}}{{i!\,\,(2m + 1 - j - i)!}}} }
\end{eqnarray*}
\emph{and}
\begin{eqnarray*}
\sum\limits_{j = 0}^{m - 1} {\frac{{( - 1)^{j + 1} }}{{(j + 1)!}}\sum\limits_{p = 1}^{2m - 1 - j} {\frac{{h^{p +
1} Z_p }}{{p!\,\,(2m - 1 - j - p)!}} } }&= & \sum\limits_{j = 2}^{m + 1} {\frac{{h^j Z_{j - 1} }}{{(j -
1)!}}\sum\limits_{i = 1}^{m} {\frac{{( - 1)^i }}{{i!\,\,(2m + 1 - j - i)!}} } }+\\
&& \sum\limits_{j = m + 2}^{2m} {\frac{{h^j Z_{j - 1} }}{{(j - 1)!}}\sum\limits_{i = 1}^{2m + 1 - j}
{\frac{{( - 1)^i }}{{i!\,\,(2m + 1 - j - i)!}}} }.
\end{eqnarray*}

The proof of Lemma 5.2 is obtained by expansion in powers of $h$
of the left sides of given identities.

For the coefficients of the optimal quadrature formulas of the
form (1.4) the following theorem holds.

\textbf{Theorem 5.4.} \textit{Among quadrature formulas of the form  (1.4) with the error functional (1.5) in
the space $L_2^{(m)}(0,1)$ for $m\geq 4$ there exists unique optimal formula which coefficients are determined by the
following formulas
$$
C[0] = h\left( {\frac{1}{2} + \sum\limits_{k = 1}^{m - 1}d_k
{\frac{{q_k^N  - q_k }}{{1 - q_k }}} } \right),\eqno (5.24)
$$
$$
C[\beta ] = h\left( {1 + \sum\limits_{k = 1}^{m - 1}d_k {(
q_k^\beta   + q_k^{N - \beta } )} } \right),\ \ \beta = \overline
{1,N - 1},\eqno   (5.25)
$$
$$
C[N] = h\left( {\frac{1}{2} + \sum\limits_{k = 1}^{m - 1}d_k
{\frac{{ q_k^N  - q_k }}{{1 - q_k }}} } \right),\eqno (5.26)
$$
$$
A = h^2 \left( {\frac{1}{{12}} - \sum\limits_{k = 1}^{m - 1}
d_k{\frac{{ q_k  + q_k^{N + 1} }}{{(1 - q_k )^2 }}} }
\right),\eqno     (5.27)
$$
$$
B = h^4 \left( {  \frac{B_4}{{4!}} - \frac{1}{{3!}}\sum\limits_{k
= 1}^{m - 1} d_k{\sum\limits_{i = 0}^{3} {\frac{{ q_k  + ( - 1)^{i + 1}q_k^{N +
i}  }}{{(q_k  - 1)^{i + 1} }}} } \Delta ^i 0^{3} }
\right),\eqno
(5.28)
$$
where $d_k$ satisfy the following system of $m - 1$ linear
equations
                $$
\sum\limits_{k = 1}^{m - 1}d_k {\sum\limits_{i = 0}^j {\frac{{q_k
+ ( - 1)^{i + 1}q_k^{N + i}  }}{{(q_k  - 1)^{i + 1} }}} } \Delta
^i 0^j  = \frac{{B_{j + 1} }}{{j + 1}},\,\,\,\,\,\,j = \overline
{4,m - 1},\eqno     (5.29)
$$
$$
\sum\limits_{k = 1}^{m - 1}d_k {\sum\limits_{i = 0}^{j}
{\frac{{q_k  + (-1)^{i + 1}q_k^{N + i} }}{{(q_k  - 1)^{i + 1} }}}
} \Delta ^i 0^{j}  = 0, \ \ \  \ j=2, 2m-4, 2m-2.\eqno       (5.30)
$$
Here $B_\alpha $ are Bernoulli numbers, $\Delta ^i \gamma ^j  $ is the finite difference of order $i$ of
$\gamma^j $, $\,\Delta ^i 0^j$ is given in Theorem 2.1, $q_k $ are given in Theorem 5.1.}

\textbf{Proof.} First we consider the first sum of  equation
(4.1). For this sum we have
\begin{eqnarray*}
 S&=& \sum\limits_{\gamma  = 0}^N {C[\gamma ]\frac{{|h\beta  - h\gamma |^{2m - 1} }}{{2(2m - 1)!}} = }\\
&=&C[0]\frac{{(h\beta )^{2m - 1} }}{{(2m - 1)!}} + \sum\limits_{\gamma  = 1}^\beta  {C[\gamma ]\frac{{(h\beta  -
h\gamma )^{2m - 1} }}{{(2m - 1)!}}}  - \sum\limits_{\gamma  = 0}^N {C[\gamma ]} \frac{{(h\beta  - h\gamma )^{2m
- 1} }}{{2(2m - 1)!}}.
\end{eqnarray*}
The last two sums of the expression $S$ we denote
$$
S_1  = \sum\limits_{\gamma  = 1}^\beta  {C[\gamma ]\frac{{(h\beta
- h\gamma )^{2m - 1} }}{{(2m - 1)!}}},\ \ \  S_2  =
\sum\limits_{\gamma  = 0}^N {C[\gamma ]} \frac{{(h\beta  - h\gamma
)^{2m - 1} }}{{2(2m - 1)!}}
$$
and we calculate them separately.\\
By using (5.17) and formulas (2.4), (2.5) for $S_1 $ we have
\begin{eqnarray*}
S_1&=&\sum\limits_{\gamma  = 0}^\beta {h\left( {1 + \sum\limits_{k = 1}^{m - 1} {\left( {d_k q_k^\gamma + p_k
q_k^{N - \gamma } } \right)} } \right)\frac{{(h\beta  - h\gamma )^{2m - 1} }}{{(2m - 1)!}} = } \\
 &=& \frac{{h^{2m} }}{{(2m - 1)!}}\left[ {\sum\limits_{\gamma  = 0}^{\beta  - 1} {\gamma ^{2m - 1}  +
\sum\limits_{k = 1}^{m - 1} {\left( {d_k q_k^\beta \sum\limits_{\gamma  = 0}^{\beta  - 1} {q_k^{ - \gamma }
\gamma ^{2m - 1} }  + p_k q_k^{N - \beta } \sum\limits_{\gamma  = 0}^{\beta  - 1} {q_k^\gamma  \gamma ^{2m - 1}
} } \right)} } } \right] =\\
 &=& \frac{{h^{2m} }}{{(2m - 1)!}}\left[ {\sum\limits_{j = 1}^{2m} {\frac{{(2m - 1)!B_{2m - j} }}{{j! \cdot (2m -
j)!}}\beta ^j }  + } \right.\sum\limits_{k = 1}^{m - 1} {\left[ {d_k q_k^\beta \left\{ {\frac{{q_k }}{{q_k -
1}}\sum\limits_{i = 0}^{2m - 1} {\frac{{\Delta ^i 0^{2m - 1} }}{{(q_k  - 1)^i }}} } \right.} \right.}  -\\
&&- \frac{{q_k^{1 - \beta } }}{{q_k  - 1}}\left. {\sum\limits_{i = 0}^{2m - 1} {\frac{{\Delta ^i \beta ^{2m - 1}
}}{{(q_k  - 1)^i }}} } \right\} + p_k q_k^{N - \beta } \left\{ {\frac{1}{{1 - q_k }}\sum\limits_{i = 0}^{2m - 1}
{\left( {\frac{{q_k }}{{q_k  - 1}}} \right)^i \Delta ^i 0^{2m - 1}  - } } \right.\\
&& - \frac{{q_k^\beta  }}{{1 - q_k }}\left. {\left. {\left. {\sum\limits_{i = 0}^{2m - 1} {\left( {\frac{{q_k
}}{{q_k  - 1}}} \right)^i \Delta ^i \beta ^{2m - 1} } } \right\}} \right]} \right].
\end{eqnarray*}
Taking into account that $q_k $ is the root of the Euler-Frobenius polynomial $E_{2m - 2}(q)$ and using formulas
(2.3), (2.6) the expression for $S_1$ we reduce to the following form
$$S_1  = \frac{{(h\beta )^{2m} }}{{(2m)!}} +
h \cdot \frac{{(h\beta )^{2m - 1} }}{{(2m - 1)!}}B_1  + h^{2m}
\sum\limits_{j = 1}^{2m - 2} {\frac{{B_{2m - j} }}{{j!(2m -
j)!}}\beta ^j  + } $$ $$ + h^{2m} \sum\limits_{j = 0}^{2m - 1}
{\frac{{\beta ^{2m - 1 - j} }}{{j!(2m - 1 - j)!}}\sum\limits_{k =
1}^{m - 1} {\sum\limits_{i = 0}^j {\frac{{ - d_k q_k  + p_k q_k^{N
+ i} ( - 1)^i }}{{(q_k  - 1)^{i + 1} }}\Delta ^i 0^j .} } }
\eqno            (5.31)
$$
Now we consider $S_2$. By using equations (4.4)-(4.6)  we rewrite
the expression $S_2$ in powers of $h\beta$
$$S_2  = \sum\limits_{\gamma  = 0}^N C[\gamma ]\frac{{(h\beta  -
h\gamma )^{2m - 1} }}{{2(2m - 1)!}} =
$$
$$
=\sum\limits_{j=0}^{m-1}\frac{(h\beta)^{2m-1-j}(-1)^{j}}
{2(2m-1-j)!j!}\sum\limits_{\gamma=0}^NC[\gamma](h\gamma)^{j}+
\sum\limits_{j=m}^{2m-1}\frac{(h\beta)^{2m-1-j}(-1)^{j}}
{2(2m-1-j)!j!}\sum\limits_{\gamma=0}^NC[\gamma](h\gamma)^{j}=
$$
$$
=\frac{(h\beta)^{2m-1}}{2(2m-1)!}-\frac{(h\beta)^{2m-2}}{4(2m-2)!}+\frac{(h\beta)^{2m-3}}{4(2m-3)!}\left(\frac{1}{3}+2A\right)-
\frac{(h\beta)^{2m-4}}{2(2m-4)!3!}\left(\frac{1}{4}+3A\right)+$$
$$+
\sum\limits_{j=4}^{m-1}\frac{(h\beta)^{2m-1-j}(-1)^{j}}{2(2m-1-j)j!}
\left(\frac{1}{j+1}+j
 A+j(j-1)(j-2)B\right)+
$$
$$
+\sum\limits_{j=m}^{2m-1}\frac{(h\beta)^{2m-1-j}(-1)^{j}}
{2(2m-1-j)!j!}\sum\limits_{\gamma=0}^NC[\gamma](h\gamma)^{j}.\eqno (5.32)
$$

Substituting (5.2) and $S$ into equation (4.1) and using (5.31),
(5.32) we have
$$\frac{{(h\beta )^{2m} }}{{(2m)!}} + C[0]\frac{{(h\beta )^{2m - 1}
}}{{(2m - 1)!}} + h\frac{{(h\beta )^{2m - 1} }}{{(2m - 1)!}}B_1  +
\sum\limits_{j = 1}^{2m - 2} {\frac{{B_{2m - j} h^{2m - j} (h\beta
)^j }}{{j!(2m - j)!}} + } $$
$$
+ \sum\limits_{j = 0}^{2m - 1} {\frac{{h^{j + 1} (h\beta )^{2m - 1
- j} }}{{j!(2m - 1 - j)!}}\sum\limits_{k = 1}^{m - 1}
{\sum\limits_{i = 0}^j {\frac{{ - d_k q_k  + p_k q_k^{N + i} ( -
1)^i }}{{(q_k  - 1)^{i + 1} }}\Delta ^i 0^j  - } } }
$$
$$
-\frac{(h\beta)^{2m-1}}{2(2m-1)!}+\frac{(h\beta)^{2m-2}}{4(2m-2)!}-\frac{(h\beta)^{2m-3}}{4(2m-3)!}\left(\frac{1}{3}+2A\right)+
\frac{(h\beta)^{2m-4}}{2(2m-4)!3!}\left(\frac{1}{4}+3A\right)-$$
$$-
\sum\limits_{j=4}^{m-1}\frac{(h\beta)^{2m-1-j}(-1)^{j}}{2(2m-1-j)j!}
\left(\frac{1}{j+1}+j
 A+j(j-1)(j-2)B\right)-
$$
$$
-\sum\limits_{j=m}^{2m-1}\frac{(h\beta)^{2m-1-j}(-1)^{j}}
{2(2m-1-j)!j!}\sum\limits_{\gamma=0}^NC[\gamma](h\gamma)^{j}
$$
$$
- A\frac{{(h\beta )^{2m -2} }}{{2(2m - 2)!}} -
 A\sum\limits_{j = 0}^{2m - 2} {\frac{{(h\beta )^{2m - 2 - j} ( -
1)^j }}{{2 \cdot j! \cdot (2m - 2 - j)!}}} - B\frac{{(h\beta )^{2m - 4} }}{{2(2m - 4)!}}
-
$$
$$
- B\sum\limits_{j = 0}^{2m - 4} {\frac{{(h\beta )^{2m - 4 - j} (-
1)^j }}{{2 \cdot j! \cdot (2m -4 - j)!}}}  + \sum\limits_{\alpha
= 0}^{m - 1} {\lambda _\alpha  (h\beta )^\alpha  }  =
\frac{{(h\beta )^{2m} }}{{(2m)!}} + \sum\limits_{j = 0}^{2m - 1}
{\frac{{( - h\beta )^{2m - 1 - j} }}{{2 \cdot (2m - 1 - j)! \cdot
(j + 1)!}}}.
$$
Hence equating coefficients of the same powers of $h\beta $ gives
$$ \sum\limits_{j= 0}^{m - 1} {\lambda _j(h\beta )^j} =\sum\limits_{j= 0}^{m - 1}\frac{(h\beta)^{j}}{j!}\left[\frac{(-1)^j}{2(2m-j)!} +\frac{B(-1)^j}{2(2m-4-j)!}+\frac{A(-1)^j}{2(2m-2-j)!} \right.-$$
$$
- \frac{h^{2m - j}}{(2m-1-j)!} \sum\limits_{k = 1}^{m - 1} {\sum\limits_{i = 0}^{2m
- 1 - j} {\frac{{ - d_k q_k  + p_k q_k^{N + i} ( - 1)^i }}{{(q_k
- 1)^{i + 1} }}\Delta ^i 0^{2m - 1 - j}  + } }
$$
$$
+\left. \frac{1}{2(2m-1-j)!}\sum\limits_{\gamma  = 0}^N {C[\gamma ]} ( -
h\gamma )^{2m - 1 - j}\right] -\sum\limits_{j = 1}^{m - 1} {\frac{{B_{2m - j} h^{2m - j} (h\beta
)^j }}{{j!(2m - j)!}}  },\eqno (5.33)
$$
$$
\sum\limits_{k = 1}^{m - 1} {\sum\limits_{i = 0}^{j+4} {\frac{{d_k q_k
+ p_k q_k^{N + i} ( - 1)^{i + 1} }}{{(q_k  - 1)^{i + 1} }}} }
\Delta ^i 0^{j+4}  = \frac{{B_{j + 5} }}{{j + 5}},\,\,\,\,\,j =
\overline {0,m - 5} ,\eqno (5.34)
$$
$$\sum\limits_{k =
1}^{m - 1} d_k{\sum\limits_{i = 0}^{2}  {{{ q_k^{N + i}( - 1)^{i +
1}  + q_k } \over {(q_k -1  )^{i + 1} }}\Delta ^i } } 0^{2} =
0,\eqno (5.35)
$$
 $$
C[0] = h\left( {\frac{1}{2} + \sum\limits_{k = 1}^{m - 1}
{\frac{{p_k q_k^N  - d_k q_k }}{{1 - q_k }}} } \right),\eqno
(5.36)
$$
$$
A = h^2 \left( {\frac{1}{{12}} - \sum\limits_{k = 1}^{m - 1}
{\frac{{d_k q_k^{}  + p_k q_k^{N + 1} }}{{(1 - q_k )^2 }}} }
\right),\eqno  (5.37)
$$
$$
B = h^4 \left( {  -\frac{1}{{720}} + \frac{1}{{3!}}\sum\limits_{k
= 1}^{m - 1} {\sum\limits_{i = 0}^{3} {\frac{{d_k q_k  + p_k
q_k^{N + i} ( - 1)^{i + 1} }}{{(q_k  - 1)^{i + 1} }}} } \Delta ^i
0^{3} } \right).\eqno (5.38)
$$
Substituting the expressions (5.36) into (4.4), also taking into
account (5.17), we find $C[N]$ which have the following form
    $$
C[N] = h\left( {\frac{1}{2} + \sum\limits_{k = 1}^{m - 1}
{\frac{{d_k q_k^N  - p_k q_k }}{{1 - q_k }}} } \right).\eqno
(5.39)
$$

Now, putting the values of $ {\lambda
_j} $  into (4.2) the following equations for we get unknowns $d_k $ and $p_k$:
$$
\sum\limits_{k = 1}^{m - 1}\sum\limits_{i = 0}^{2m-2} {\frac{{d_kq_k^{i} + p_kq_k^{N + 1} ( -
1)^{i + 1} }}{{(1-q_k)^{i + 1} }}\Delta ^i 0^{2m-2}  }  =\sum\limits_{k = 1}^{m - 1} \sum\limits_{i = 0}^{2m-2}  {\frac{{d_kq_k^{N+i}  +
p_kq_k(-1)^{i + 1} }}{{(1 - q_k)^{i + 1} }}\Delta ^i
0^{2m-2}  },\eqno (5.40)
$$
$$\sum\limits_{k = 1}^{m - 1} \sum\limits_{i = 0}^{j-1}  {\frac{{d_kq_k^{N+i}  +
p_kq_k(-1)^{i + 1} }}{{(1 - q_k)^{i + 1} }}\Delta ^i
0^{j-1}  }=\frac{B_j}{j},\eqno (5.41)
$$
$$
\sum\limits_{k = 1}^{m - 1}\sum\limits_{i = 0}^{3} {\frac{{d_kq_k + p_kq_k^{N + i} ( -
1)^{i + 1} }}{{(q_k-1)^{i + 1} }}\Delta ^i 0^{3}  }  =\sum\limits_{k = 1}^{m - 1} \sum\limits_{i = 0}^{3}  {\frac{{d_kq_k^{N+i}  +
p_kq_k(-1)^{i + 1} }}{{(1 - q_k)^{i + 1} }}\Delta ^i
0^{3}  },\eqno (5.42)
$$
    $$
\sum\limits_{k = 1}^{m - 1} {\sum\limits_{i = 0}^{2}
{\frac{{d_k q_k^{N+i}  + p_k q_k( - 1)^{i + 1} }}{{(1-q_k)^{i + 1} }}} } \Delta ^i 0^{2}  = 0,\eqno   (5.43)
$$
$$
\sum\limits_{k = 1}^{m - 1}{\frac{{d_kq_k + p_kq_k^{N + 1}  }}{{(q_k-1)^{2} }}}  =\sum\limits_{k = 1}^{m - 1} {\frac{{d_kq_k^{N+1}  +
p_kq_k}}{{(1 - q_k)^{2} }}},\eqno (5.44)
$$
$$
\sum\limits_{k = 1}^{m - 1}\sum\limits_{i = 0}^{2m-4} {\frac{{d_kq_k^{i} + p_kq_k^{N + 1} ( -
1)^{i + 1} }}{{(1-q_k)^{i + 1} }}\Delta ^i 0^{2m-4}  }  =\sum\limits_{k = 1}^{m - 1} \sum\limits_{i = 0}^{2m-4}  {\frac{{d_kq_k^{N+i}  +
p_kq_k(-1)^{i + 1} }}{{(1 - q_k)^{i + 1} }}\Delta ^i
0^{2m-4}  },\eqno (5.45)
$$
thus, from equalitions (5.34)- (5.35) and (5.40)- (5.45) we get
$$
\sum\limits_{k = 1}^{m - 1}(d_k-p_k) {\sum\limits_{i = 0}^{2m-2} {\frac{{ q_k
+ q_k^{N + i} ( - 1)^{i + 1} }}{{(q_k  - 1)^{i + 1} }}} }
\Delta ^i 0^{2m-2}  = 0 ,\eqno (5.46)
$$
$$
\sum\limits_{k = 1}^{m - 1}(d_k-p_k) {\sum\limits_{i = 0}^{2m-4} {\frac{{ q_k
+ q_k^{N + i} ( - 1)^{i + 1} }}{{(q_k  - 1)^{i + 1} }}} }
\Delta ^i 0^{2m-4}  = 0 ,\eqno (5.47)
$$
$$
\sum\limits_{k = 1}^{m - 1}(d_k-p_k) {\sum\limits_{i = 0}^{2} {\frac{{ q_k
+ q_k^{N + i} ( - 1)^{i + 1} }}{{(q_k  - 1)^{i + 1} }}} }
\Delta ^i 0^{2}  = 0 ,\eqno (5.48)
$$
$$
\sum\limits_{k = 1}^{m - 1}(d_k-p_k) {\sum\limits_{i = 0}^{j} {\frac{{q_k
+  q_k^{N + i} ( - 1)^{i + 1} }}{{(q_k  - 1)^{i + 1} }}} }
\Delta ^i 0^{j}  = \frac{{B_{j + 1} }}{{j + 1}},\,\,\,\,\,j =
\overline {4,m - 1} ,\eqno (5.49)
$$

Taking into account uniqueness of the optimal coefficients, we
conclude, that the homogeneous system of linear equations
(5.46)-(5.49) has trivial solution. This
means, that
$$
d_k=p_k,\ \ \ k=1,2,...,m-1.\eqno (5.50)
$$
Then, using (5.50), from (5.46)- (5.49) we get (5.29), (5.30), and
from (5.17), (5.36)-(5.39) we obtain (5.24)-(5.28).

Theorem 5.4 is proved.

\subsection{The norm of the error functional of optimal quadrature
formulas of the form (1.4)}

For square of the norm of the error functional (1.5) of optimal
quadrature formulas of the form (1.4) the following holds

\textbf{Theorem 5.5.} \textit{For square of the norm of the error
functional  (1.5) of the optimal quadrature formula of the form
(1.4) on the space $L_2^{(m)}(0,1)$ for $m\geq 4$ the following holds
    $$
\left\|\stackrel{\circ}{\ell}|L_2^{(m)*}(0,1)\right\|^2  = ( -
1)^{m + 1} \left[ {\frac{{B_{2m} h^{2m} }}{{(2m)!}}} \right. -
\frac{2{h^{2m + 1} }}{{(2m)!}}\left. {\sum\limits_{k = 1}^{m -
1}d_k {\sum\limits_{i = 0}^{2m} {\frac{{q_k +q_k^{N+i} ( -
1)^{i+1} }}{{(q_k -1)^{i + 1} }}} } \Delta ^i 0^{2m} } \right],
$$
where $d_k$ are determined from the system (5.29)-(5.30), $B_{2m}
$ are Bernoulli numbers, $\Delta^i0^{2m}$ is given in Theorem 2.1,
$q_k $ are given in Theorem 5.1. }

\textbf{Proof.} Computing definite integrals in the expression
(3.5) of $||\ell||^2 $ we get
$$
\left\| {\ell } \right\|^2  = ( - 1)^{m } \left[
{\sum\limits_{\beta  = 0}^N {C [\beta ]\left(
{\sum\limits_{\gamma  = 0}^N {C [\gamma ]} } \right.} }
\right.{{|h\beta  - h\gamma|^{2m - 1} }
\over {2(2m - 1)!}} -
$$
$$
\left.- 2\int\limits_0^1 {{{|x - h\beta |^{2m - 1} }
\over {2(2m - 1)!}}} dx-2AP_{2m-2}(h\beta)-
 2BP_{2m-4}(h\beta)\right)+
$$
$$
+\frac{2A}{(2m-1)!}+\frac{2B+A^2}{(2m-3)!}+
 \left.
 \frac{2AB}{(2m-5)!}+\frac{B^2}{(2m-7)!}+\frac{1}{(2m+1)!}
\right].
$$
where $f(h\beta)$ is defined by formula (5.2).  As is obvious from
here according to (4.1) the expression into curly brackets is
equal to the polynomial $ \sum\limits_{\alpha =
0}^{m - 1} {\lambda _\alpha (h\beta )^\alpha } $. Then $||\ell
||^2 $ has the form
$$
\left\| {\ell } \right\|^2  = ( - 1)^{m } \left[
{\sum\limits_{\beta  = 0}^N {C [\beta ]\left(
- \sum\limits_{j =
0}^{m - 1} {\lambda _j (h\beta )^j }-f(h\beta) \right)} }\right.
 -
$$
$$
-A\sum\limits_{\beta  = 0}^N C [\beta ]P_{2m-2}(h\beta)-
 B\sum\limits_{\beta  = 0}^N C [\beta ]P_{2m-4}(h\beta)+
$$
$$
+\frac{2A}{(2m-1)!}+\frac{2B+A^2}{(2m-3)!}+
 \left.
 \frac{2AB}{(2m-5)!}+\frac{B^2}{(2m-7)!}+\frac{1}{(2m+1)!}
\right].
$$
Hence using (4.2) and (4.3) we get
$$
\left\| {\ell } \right\|^2  = ( - 1)^{m } \left[
- \sum\limits_{j =
0}^{m - 1} {\lambda _j \sum\limits_{\beta  = 0}^N C [\beta ](h\beta )^j }-\sum\limits_{\beta  = 0}^N C [\beta ]f(h\beta)  \right.
 +
$$
$$
+A\sum\limits_{j =
2}^{m - 1} j\lambda _j+B\sum\limits_{j =
4}^{m - 1}j(j-1)(j-2)\lambda _j+
$$
$$
+\frac{A}{(2m-1)!}+\frac{B}{(2m-3)!}+
 \left.
 \frac{1}{(2m+1)!}
\right]. \eqno   (5.51)
$$

From (5.51) after some simplifications, using (5.33), (5.2) and (4.4)-(4.6)
we
have
$$\left\| {\ell } \right\|^2  = ( -
1)^{m }\left[ -\sum\limits_{j = 0}^{m - 1} {(-1)^{j+1} \over
(j+1)!(2m-1-j)!}\sum\limits_{\gamma  = 0}^N {C [\gamma ]( h\gamma
)^{2m - 1 - j}   } -\right.
$$
$$-\sum\limits_{j = 0}^{m - 1}\frac{h^{2m
-j}}{(j+1)!(2m-1-j)!}{ \sum\limits_{k = 1}^{m - 1}d_k
{\sum\limits_{i = 0}^{2m - 1 - j} {{{ q_k  + q_k^{N + i } ( -
1)^{i+ 1} } \over {(q_k  - 1)^{i + 1} }}} \Delta ^i 0^{2m - 1 - j}
 } }
$$
$$
- \sum\limits_{j = 0}^{m - 1} \frac{(-1)^j}{(j+1)!(2m-j)!}-
\sum\limits_{j = 0}^{m - 5} \frac{B(-1)^j}{(j+1)!(2m-4-j)!}-$$
$$-
\sum\limits_{j = 0}^{m - 3} \frac{A(-1)^j}{(j+1)!(2m-2-j)!}
+\sum\limits_{j = 1}^{m - 1}
\frac{B_{2m-j}h^{2m-j}}{(j+1)!(2m-j)!}
$$
$$\left.
-\frac{1}{(2m)!}\sum\limits_{\beta  = 0}^N C [\beta
](h\beta]^{2m}+\frac{A}{(2m-1)!}+\frac{B}{(2m-3)!} +
\frac{1}{(2m+1)!}\right] \eqno   (5.52)$$

When $\alpha>m - 1$ using (5.17) and formulas (2.4)-(2.6) we get
$$
\sum\limits_{\gamma  = 0}^N {C[\gamma ](h\gamma )^\alpha   =
\frac{1}{{\alpha  + 1}} + \sum\limits_{j = 1}^{\alpha  - 1}
{\frac{{\alpha !B_{\alpha  + 1 - j} }}{{j!(\alpha  + 1 -
j)!}}h^{\alpha  + 1 - j}  + } }
$$
$$
+ h^{\alpha  + 1} \sum\limits_{k = 1}^{m - 1} {\sum\limits_{i =
0}^\alpha  {\frac{{d_k q_k^i  + p_k q_k^{N + 1} ( - 1)^{i + 1}
}}{{(1 - q_k )^{i + 1} }}\Delta ^i 0^\alpha   - } }
$$
$$
- \sum\limits_{j = 1}^\alpha  {\frac{{\alpha !h^{j + 1}
}}{{j!(\alpha  - j)!}}\sum\limits_{k = 1}^{m - 1} {\sum\limits_{i
= 0}^j {\frac{{d_k q_k^{N + i}  + p_k q_k^{} ( - 1)^{i + 1} }}{{(1
- q_k )^{i + 1} }}\Delta ^i 0^j } } }.\eqno   (5.53) $$ Using
Lemmas 5.1, 5.2 and taking into account (2.1), (2.3), after some
simplifications, from (5.52) we have

$$\left\| {\ell } \right\|^2  = ( -
1)^{m }\left[K_1 +K_2-\frac{B_{2m}h^{2m}}{(2m)!}+\frac{h^{2m+1}Z_{2m}}{(2m)!}\right.
$$
$$-\frac{h^{2m+1}}{(2m)!}\sum\limits_{k = 1}^{m - 1}d_k {\sum\limits_{i =
0}^{2m}  {\frac{{ q_k^i  + q_k^{N + 1} ( - 1)^{i + 1}
}}{{(1 - q_k )^{i + 1} }}\Delta ^i 0^{2m}  } }
-
$$
$$-\sum\limits_{j = 2}^{m} {B_jh^j \over
j!}\sum\limits_{i  = 0}^m \frac{(-1)^i}{i!(2m+1-i-j)!}
+ \sum\limits_{j =2}^{m} {\frac{{h^jZ_{j-1}
}}{{(j-1) !}}\sum\limits_{i = 0}^{m} {\frac{{( - 1)^i
}}{{i!(2m+1-i-j)!}}} }
$$
$$\left.+
\sum\limits_{j = 0}^{m - 4} \frac{B(-1)^j}{j!(2m-3-j)!}+
\sum\limits_{j = 0}^{m - 2} \frac{A(-1)^j}{j!(2m-1-j)!}
\right] \eqno   (5.54)
$$

where
$$ K_1=\sum\limits_{j = m+1}^{2m}\frac{h^jZ_{j-1}}{(j-1)!}\sum\limits_{i = 1}^{2m+1-j}\frac{(-1)^i
}{i!(2m+1-i-j)!}+\sum\limits_{j = m+1}^{2m}\frac{h^jZ_{j-1}}{(j-1)!(2m+1-j)!}=$$
$$
=\sum\limits_{j = m+1}^{2m}\frac{h^jZ_{j-1}}{(j-1)!}\sum\limits_{i = 0}^{2m+1-j}\frac{(-1)^i
}{i!(2m+1-i-j)!}=\sum\limits_{j = m+1}^{2m}\frac{h^jZ_{j-1}}{(j-1)!}\frac{(1-1)^{2m+1-j}}{(2m+1-j)!}=0,
$$
and
$$
K_2=-\sum\limits_{j = m+1}^{2m-1}\frac{h^jB_{j}}{j!}\sum\limits_{i = 0}^{2m-j}\frac{(-1)^i
}{i!(2m+1-i-j)!}+\sum\limits_{j = m+1}^{2m-1}\frac{h^jB_{j}}{j!(2m+1-j)!}=
$$
$$
-\sum\limits_{j = m+1}^{2m-1}\frac{h^jB_{j}}{j!}\left(\sum\limits_{i = 0}^{2m-j}\frac{(-1)^i
}{i!(2m+1-i-j)!}-\frac{1}{(2m+1-j)!}\right)=
$$
$$
=-\sum\limits_{j = m+1}^{2m-1}\frac{h^jB_{j}}{j!}\left(\frac{(-1)^j-1}{(2m+1-j)!}\right)=0
$$
Since $K_1= K_2=0$ then from (5.54) we get the folloving

$$\left\| {\ell } \right\|^2  = ( -
1)^{m }\left[-\frac{B_{2m}h^{2m}}{(2m)!}+\frac{{h^{2m + 1} }}{{(2m)!}} {\sum\limits_{k = 1}^{m -
1}d_k {\sum\limits_{i = 0}^{2m} {\frac{{q_k +q_k^{N+i} ( -
1)^{i+1} }}{{(q_k -1)^{i + 1} }}} } \Delta ^i 0^{2m} } \right.
$$
$$\left. -\frac{h^{2m+1}}{(2m)!}\sum\limits_{k = 1}^{m - 1}d_k {\sum\limits_{i =
0}^{2m}  {\frac{{ q_k^i  + q_k^{N + 1} ( - 1)^{i + 1}
}}{{(1 - q_k )^{i + 1} }}\Delta ^i 0^{2m}  } }\right]
$$

Hance using Lemma 5.1 we have the statimant of the theorem Theorem 5.5 is proved.

\textbf{Corollary 5.1.} \textit{In the space $L_2^{(4)}(0,1)$
among quadrature formulas of the form (1.4) with the error
functional (1.5) there exists unique optimal formula  whose
coefficients are determined by the following formulas}
    $$
C[\beta] = \left\{
\begin{array}{ll}
\frac{h}{2},& \beta=0,N,\\
h,&\beta=\overline{1,N-1},\\
\end{array} \right.
$$
$$
A =\frac{h^2}{12},\ \ \  B =-\frac{h^4}{720}.
$$
\textit{Furthermore for square of the norm of the error functional
the following is valid}
$$
\left\|\stackrel{\circ}{\ell}|L_2^{(4)*}(0,1) \right\|^2  =
\frac{h^8}{1209600 }.
$$

\textbf{Corollary 5.2.} \textit{In the space $L_2^{(5)}(0,1)$
among quadrature formulas of the form (1.4) with the error
functional (1.5) there exists unique optimal formula  whose
coefficients are determined by the following formulas }
    $$
C[\beta] = \left\{
\begin{array}{ll}
\frac{h}{2},& \beta=0,N,\\
h,&\beta=\overline{1,N-1},\\
\end{array} \right.
$$
$$
A =\frac{h^2}{12},\ \ \  B =-\frac{h^4}{720}.
$$
\textit{Furthermore for square of the norm of the error functional the following is valid}
$$
\left\|\stackrel{\circ}{\ell}|L_2^{(5)*}(0,1) \right\|^2  =
\frac{h^{10}}{47900160}.
$$

\section*{Acknowledgements}
The work of the third author was supported in part by the grant  YoF4-OT-010509-YoF4-5 for young scientists.

\end{document}